\documentclass[12pt]{amsart}
\usepackage{amsmath,amscd,amssymb, amsbsy,eufrak,euscript}
\usepackage{syntonly}

\def\al{\alpha}

\def\be{\beta}
\def\ep{\varepsilon}
\def\de{{\delta}}

\def\la{\lambda}

\def\ga{\gamma}

\def\si{\sigma}

\def\de{{\delta}}

\def\cD{\mathcal D}
\def\cE{\mathcal E}
\def\cG{\mathcal G}

\def\cV{\mathcal V}

\def\cS{\mathcal S}

\def\limn{\lim_{n\to\infty}}

\def\l{\ell}
\newcommand{\comment}[1]{}
\newcommand{\bq}{\begin{equation}}
\newcommand{\tq}{\end{equation}}

\newcommand{\bs}{\begin{split}}
\newcommand{\es}{\end{split}}

\def\al{\alpha}
\def\be{\beta}
\def\ga{\gamma}

\def\de{\delta}

\def\ep{\varepsilon}

\def\l{\ell}
\def\la{\lambda}
\def\si{\sigma}

\def\limn{\lim_{n \to \infty}}

\def \cE{\mathcal E}
\def\cS{\mathcal S}
\numberwithin{equation}{section}

\newtheorem{thm}{Theorem}
\newtheorem{lm}{Lemma}

\theoremstyle{remark}

\newtheorem{con}{}
\def\lbl(#1){}
\def\EQN(#1){\lbl(#1)
\protect\begin{equation}\protect\label{#1}}

\def\CON(#1){\lbl(#1) \begin{con}
\label{#1} \end{con}}

\title{On the range of the simple random walk bridge on groups}
\author{Itai Benjamini \and  Roey Izkovsky \and Harry Kesten}
\date{Jan 2006}
\address{The Weizmann Institute of Science\\76100 Rehovot\\Israel}
\email{itai.benjamini@gmail.com, roey.izkovsky@gmail.com}
\address{Department of Mathematics\\
Malott Hall\\
Cornell University\\
Ithaca NY 14853\\USA} \email{kesten@math.cornell.edu}

\begin{document}
\begin{abstract}  Let $\cG$ be a vertex transitive graph. A study of the range
of simple random walk on $\cG$ and of its bridge is proposed.
While it is expected that on a graph of polynomial growth the
sizes of the range of the unrestricted random walk and of its
bridge are the same in first order, this is not the case on some
larger graphs such as regular trees. Of particular interest is the
case when $\cG$ is the Cayley graph of a group $G$. In this case
we even study the range of a general symmetric (not necessarily
simple) random walk on $G$. We hope that the few examples for
which we calculate the first order behavior of the range here will
help to discover some relation between the group structure and the
behavior of the range. Further problems regarding bridges are
presented.
\end{abstract}

\maketitle
\pagestyle{myheadings}
\markboth{Itai Benjamini, Roey
Izkovsky and Harry Kesten} {On the range of the simple random
walk bridge on groups}

\section{{\bf Introduction.}}

A simple random walk bridge of length $n$ on a graph, is a simple
random walk (SRW) conditioned to return to the starting point of the
walk at time $n$. In this note we initiate a study of bridges on
vertex transitive graphs, concentrating mainly on the range of a
bridge. There is a considerable literature (see for instance
\cite{DvE},
\cite{Sp1}, \cite{Sp2}, \cite{JP}, \cite{De}, \cite{DoV},
\cite{Ha}) on the range of a random walk on $\mathbb Z^d$ and on
more general graphs. The first result in this area seems to be the
following strong law of large numbers from \cite{DvE},
 \cite{Sp1}, Theorem 4.1 : Let $\{S_n\}_{n \ge 0}$ be a random walk on
$\mathbb Z$ and let $R_n :=\big|\{S_0, S_1,\dots, S_{n-1}\}\big|$
be its range at time $n$. Then
\begin{equation}
\limn \frac 1n R_n \to 1-F \text{ a.s.}, \label{1.1}
\end{equation}
where
\begin{equation}
F := P\{S_n = S_0 \text{ for some } n \ge 1\}. \label{1.2}
\end{equation}
This result is for an unrestricted random walk, that is, for $S_n =
\sum_{i=1}^n X_i$ with the $X_i$ i.i.d. $\mathbb Z$-valued random
variables. It was extended in \cite{Sp2}, \cite{De} to the case when
the $\{X_i\}$ form a stationary ergodic sequence. The proof is a
simple application of Kingman's subadditive ergodic theorem. It can
even be extended to a simple random walk on a vertex transitive graph
(see below for a definition). In this paper we are interested in
comparing the limit in (\ref{1.1}), (\ref{1.2}) with the limit of
$(1/n)R_n$ when $\{S_0, \dots, S_n\}$ is conditioned on the event
$\cE_n :=\{S_n = S_0\}$. In this conditioned case, which has the
condition varying with $n$, we can only speak of the limit in
probability, since an almost sure limit is meaningless. In a number
of examples we shall calculate this limit in probability of
$(1/n)R_n$ and see that it equals $1-F$ in some cases and differs
from $1-F$ in other cases. ``Usually'' the limit of $(1/n)R_n$ under
the condition $\cE_n$ is less than or equal to the limit for
unrestricted random walk. The idea is that conditioning on $\cE_n$
will pull in $S_i$ closer to its starting point than in the
unconditioned case, and that this may diminish the range. We shall be
particularly interested in the case when $\cG$ is the Cayley graph of
a finitely generated, infinite group. One would hope that in this
case the values of the different limits for $(1/n)R_n$ give some
information about the size or structure of the group. Even though it
is unclear to what extent such group properties influence $R_n$, it
is likely that the volume growth of the group play a role (see also
Open Problem 1 later in this section). It will also be apparent from
our calculations that the behavior of the Green function $P\{S_n =
S_0\}= P\{\cE_n\}$ is significant.

Various other papers have discussed bridges of random walks on
graphs and in particular Cayley graphs. \cite{NGV} and \cite{BJ}
prove invariance principles for such bridges. \cite{E} discusses
the graph distance between the starting point of a bridge and its
``midpoint'' (to be more specific, if the bridge returns to its
starting point at time {2n}, then by its midpoint we mean the
position of the bridge at time $n$); see also the discussion
preceding and following \eqref{1.18k} below). \cite{Y} studies
still other aspects of bridges of random walks on Cayley graphs
and their relation to group structure. In particular, this
reference considers the expected value of the so-called Dehn's
function of a bridge. We shall mention some further aspects of the
range and bridges, as well as some open problems towards the end
of this introduction. In fact, some of those remarks served as
motivations for the present study.

Here is a formal description of
our set up. A countable graph $\cG$ is {\it vertex transitive} if for
any two of its vertices $v'$ and $v''$, there is a graph automorphism
$\Phi(\cdot)=
\Phi(\cdot;v',v'')$ which maps $v'$ to $v''$. Throughout we let
$\cG$ be a countably infinite, connected vertex transitive graph, all
of whose vertices have degree $\cD < \infty$ and let $e$ be a specific
(but arbitrary) vertex of $\cG$. {\it Simple random walk} on $\cG$ is
the Markov chain $\{S_n\}_{n \ge 0}$ which moves from a vertex $v$ to
any one of the neighbors of $v$ with probability $1/\cD$. More
formally, its transition probabilities are $P\{S_{n+1} = w|S_n = v\}
= 1/\cD$ if $w$ is a neighbor of $v$, and 0 otherwise. Unless stated
otherwise, we assume that $S_0 = e$. Of particular interest is the
case when $\cG$ is the {\it Cayley graph} of a finitely generated
infinite group $G$. Let $G$ be generated by the finite set
$\cS=\{g_1,\cdots,g_s,g_1^{-1}, \dots, g_s^{-1}\}$ of its elements
and their inverses. We can then take for $\cG$ the graph whose
vertices are the elements of $G$ and with an edge between $v'$ and
$v''$ if and only if $v'' = v'g_i$ or $v'' = v' g_i^{-1}$ for some $1
\le i\le s$. This graph $\cG$ depends on $\cS$ and it will be denoted by
$(G,\cS)$. $\cG$ is called the Cayley graph of $G$ corresponding to
the generating set $\cS$.

If $\al_i \ge 0, \sum _{i=1}^{2s} \al_i= 1$, we can define a random
walk $\{S_n\}_{n \ge 0}$ as follows: Let $X, X_1, X_2,
\dots$ be i.i.d. $G$-valued random variables with the distribution
\begin{equation}
P\{X= g_i\} = \al_i, \quad P\{X=g_i^{-1}\} = \al_{s+i}, \quad 1
\le i \le s.
\end{equation}
Take  $e$ to be the identity element of $G$ and set $S_n = X_1X_2
\cdots X_n$. This so-called right random walk on $G$ has transition
probabilities
\begin{equation}
P\{S_{n+1} = w\big|S_n = v\} = \sum_{i:\;vg_i = w} \al_i +
\sum_{i:\;vg_i^{-1} = w} \al_{s+i}.
\end{equation}
We shall restrict ourselves here to the symmetric case in which
\begin{equation}
\al_i = \al_{s+i} \text{ or }P\{X = g_i\} = P\{X = g_i^{-1}\}.
\label{1.5}
\end{equation}
We shall further assume that
\begin{equation}
\al_i > 0 \text{ for all $i$ and $\cS$ generates }G. \label{1.6}
\end{equation}
(Note that this condition is harmless. If it does not hold from the
start we can simply replace $G$ by the group generated by the $g_i$
and $g_i^{-1}$ with $\al_i > 0$.)

Throughout we use the following notation (this does not require
$\cG$ to be a Cayley graph): $\cE_n = \{S_n = e\}$,
$$
u_n = P\{S_n = e\} = P\{\cE_n\},
$$
\begin{equation}
f_n = P\{S_k \ne e, 1 \le k \le n-1, S_n = e\},\; F =
\sum_{n=1}^\infty f_n. \label{1.5a}
\end{equation}
$f_n$ is the probability that $S_.$ returns to $e$ for the first
time at time $n$, and $F$ is the probability that $S_.$ ever
returns to $e$. Finally, $R_n = \big|\{S_0, S_1, \dots
S_{n-1}\}\big|$.

A minor nuisance is possible periodicity of the random walk. The
{\it period} is defined as
\begin{equation}
p = \text{g.c.d.}\{n: u_n > 0\}.
\end{equation}
Since the random walk can move from a vertex $v$ to a neighbor $w$
at one step and then go back in the next step from $w$ to $v$ with
positive probability, we always have $u_2 > 0$. Thus the period is
either 1 or 2. In the latter case we have by definition
$P\{\cE_n\}= u_n = 0$ for all odd $n$. In this case it makes
little sense to talk about conditioning on the occurrence of
$\cE_n$ for odd $n$. {\it If the period is 2 all statements which
involve conditioning on $\cE_n$ shall be restricted to even $n$.}

Our first result states that under the mild condition that $u_n$ does
not tend to 0 exponentially fast (see \eqref{1.10}), $R_n$
conditioned on $\cE_n$, is in some sense no bigger than $R_n$ without
the conditioning. In the second theorem we give
sufficient conditions for the limit in probability of $(1/pn)R_{pn}$,
conditioned on $\cE_{pn}$, to equal $1-F$, which is the same as the
almost sure limit of $(1/n)R_n$ without any conditioning
(recall \eqref{1.1}). The last
theorem gives another set of sufficient conditions for the existence
of the limit in probability of $(1/pn)R_{pn}$, conditioned on
$\cE_{pn}$. However, under the conditions of Theorem 3 this limit
will often differ from $1-F$. Examples of random walks satisfying the
conditions of Theorems 2 and 3 are given after the theorems.
\begin{thm} Assume that $\{S_n\}$ is simple random walk on a vertex
  transitive graph $\cG$ or a random walk on a Cayley graph for which
\eqref{1.5} and \eqref{1.6} hold. Assume further that
\begin{equation}
\limsup_{n \to \infty}\; [u_{2n}]^{1/n} = 1. \label{1.10}
\end{equation}
Then for all $\ep > 0$
\begin{equation}
\lim_{n \to \infty, p|n} P\{\frac 1nR_n > 1-F +\ep\big|\cE_n\} =
0. \label{1.11}
\end{equation}
In particular, if $\{S_n\}$ is recurrent, then $(1/n)R_n$
conditioned on $\cE_n$ tends to 0 in probability as $n \to \infty$
through multiples of the period $p$.
\end{thm}

\begin{thm}
Assume that $\{S_n\}$ is a random walk on an infinite Cayley graph
for which \eqref{1.5} and \eqref{1.6} hold. Assume further that
\begin{equation}
\limsup_{n \to \infty} \frac{u_{2n}}{u_{4n}} < \infty.
\label{1.16a}
\end{equation}
Then for all $\ep > 0$
\begin{equation}
\lim_{n \to \infty, p|n} P\{\big|\frac 1n R_n - 1+F\big| >
\ep\big|\cE_n\} = 0. \label{1.16}
\end{equation}
\eqref{1.16} is also valid if there exist two functions $g,h \ge
0$ on $\Bbb Z_+$ which satisfy
\begin{equation}
g(n) \text{ is nondecreasing and tends to }\infty, \label{1.12}
\end{equation}
\begin{equation}
\limn \frac 1n g(n) = 0, \label{1.13}
\end{equation}
and
\begin{equation}
\limn nh\big(\big\lfloor \frac n{g(n/2)}\big\rfloor\big) = 0, \label{1.14}
\end{equation}
and are such that for all large $n$
\begin{equation}
e^{-g(n)} \le u_{2n} \le h(n). \label{1.15}
\end{equation}
\end{thm}
Note that we do not require \eqref{1.16a} in  case when
\eqref{1.12}-\eqref{1.15} hold.

\medskip
\noindent {\bf Examples.}
\newline
(i) Let $\{S_n\}$ be a random walk on a Cayley graph $(G, \cS)$
  for which \eqref{1.5} and \eqref{1.6} hold. If $G$ has polynomial
  growth, then \eqref{1.16} holds. To specialize even further,
  \eqref{1.16} holds for simple random walk on $\Bbb Z^d$. To show
  this we apply Theorem 5.1 of \cite{HeS}. This
  tells us that if $G$ has polynomial growth of order $D$, then
\begin{equation}
u_{2n} \asymp n^{-D/2},
\label{1.18}
\end{equation}
where $a(n) \asymp b(n)$ for positive $a(\cdot), b(\cdot)$ means that
there exist constants $0 < C_1 \le C_2 < \infty$ such that $C_1a(n)
\le b(n) \le C_2b(n)$ for large $n$. \eqref{1.18} trivially implies
\eqref{1.16a} and hence \eqref{1.16} (by Theorem 2). We point out
that for random walks on a Cayley graph $(G,\cS)$ which satisfy
\eqref{1.5} and \eqref{1.6}, \eqref{1.16a} is actually equivalent to
polynomial growth of $G$, or more precisely, polynomial growth of the
volume function
\begin{equation}
\begin{split}
\cV(n) = \;& \cV(n;G,\cS):= \text{number of elements of $G$ which can
be written}\\
&\text{as $h_1\cdot h_2 \cdots h_k$ with $k \le n$ and each $h_i \in \cS$ or
  $h_i^{-1} \in \cS$}.
\label{2.89ab}
\end{split}
\end{equation}
(see Lemma 4 in the next section for a proof).

\smallskip
\noindent
(ii) As we saw at the end of the preceding example, Theorem 2 deals
with random walks on Cayley graphs $(G,\cS)$ in which $G$ has
polynomial growth. As we shall see, Theorem 3 deals with some cases
in which $G$ has exponential growth. It is therefore of interest to
also look at groups of so-called intermediate groups, as constructed
by Grigorchuk in \cite{Gr}. These are finitely generated groups for
which there exist constants $0 < \al \le \be < 1$ and constants $0<
C_3,C_4 < \infty$ such that
\begin{equation}
C_3e^{n^\al} \le \cV(n) \le C_4e^{n^\be}, \quad n \ge 1.
\label{1.18cd}
\end{equation}
A random walk on a Cayley graph $(G,\cS)$ for such a group $G$ and
satisfying \eqref{1.5} and \eqref{1.6} will have
\[
u_{2n} \le C_5 \exp[-C_6 n^{\al/(\al +2)}],
\]
by virtue of Theorem 4.1 in \cite{HeS}. Moreover, since $S_n$ is always
a product of at most $n$ elements of $\cS$ or inverses of such
factors, it holds for some
$v_n \in G$ that $P\{S_n = v_n\} \ge 1/\cV(n)$ and
consequently
\[
u_{2n} \ge P\{S_n = v_n\}P\{S_n = v_n^{-1}\} \ge [\cV(n)]^{-2}
\ge C_4^{-2}\exp[-2n^\be].
\]
Thus Theorem 2 applies with the choices
$g(n) = 2n^\be +2\log C_4$
and $h(n) = C_5 \exp[-C_6 n^{\al/(\al+2)}]$.
Accordingly, \eqref{1.16} holds
for such random walks.

\smallskip
\noindent (iii) Let $G$ be a wreath product $K \wr \Bbb Z^d$  with
$K$ a finitely generated group of polynomial growth and $\{S_n\}$
a random walk on a Cayley graph $(G,\cS)$ for which \eqref{1.5}
and \eqref{1.6} hold. A specific case is the traditional
lamplighter group $\Bbb Z_2 \wr \Bbb Z$ (see \cite{PS})
for a definition of a wreath product). Then
\eqref{1.16} holds. Indeed, Theorem 3.11 in \cite{PS}
shows that there exist constants $0 < C_i < \infty$
such that
\begin{equation}
C_7\exp[-C_8n^{1/3}(\log n)^{2/3}] \le u_{2n}
\le C_9\exp[-C_{10}n^{1/3}(\log n)^{2/3}].
\label{1.18ab}
\end{equation}
Thus, \eqref{1.12}-\eqref{1.15} hold with $g(n) = C_8n^{1/3}(\log
n)^{2/3} - \log C_7$ and $h(n) = C_9\exp[-C_{10}n^{1/3}(\log n)^{2/3}]$. Note
that this argument also works if $G = \Bbb Z \wr \Bbb Z$ or for $G
= K \wr Z^d$ with $K$ finite. In the latter case we have to use
Theorem 3.5 in \cite{PS} instead of Theorem
3.11. The Remark on p. 968 of \cite{PS} leads
to many more examples to which Theorem 2 applies.

\medskip
The relation \eqref{1.16} is no longer true for simple random walk
on a regular tree (which includes the case of random walk on the
Cayley graph of a free group) as the following theorem shows.
\begin{thm}  Assume that $\{S_n\}$ is simple random walk on a vertex
  transitive graph $\cG$ or a random walk on a Cayley graph for which
\eqref{1.5} and \eqref{1.6} hold. Let $\rho$ be the radius of
convergence of the powerseries
  $U(z) := \sum_{n=0}^\infty u_nz^n$ and let
$F(z)= \sum_{n=1}^\infty f_nz^n$. Then
\begin{equation}
1 \le \rho < \infty,\; \limn [u_{pn}]^{1/(pn)} = \frac 1\rho
\text{ and } F(\rho) \le 1.
\label{1.18c}
\end{equation}
If
\begin{equation}
\text{for all $0 < \eta < 1$, } \limsup_{n \to \infty, p|n}\;
\sup_{1 \le r < (1-\eta)n, p|r} \frac {u_{n-r}}{\rho^r u_n} < \infty,
\label{1.18e}
\end{equation}
and
\begin{equation}
\text{for each fixed $r$ with $p|r$, }\lim_{n \to \infty, p|n}
\frac {u_{n-r}}{\rho^ru_n} = 1, \label{1.18ea}
\end{equation}
then for all $\ep > 0$
\begin{equation}
\limn P\{\big|\frac 1{2n} R_{2n} - \big(1-F(\rho)\big) \big|
> \ep\big|\cE_n\} = 0.
\label{1.18b}
\end{equation}
\end{thm}
\noindent {\bf Examples.} Let $\cG_b$ be a regular $b$-ary tree
(in which each vertex has degree $b+1$). Then, if $\{S_n\}$ is
simple random walk on $\cG_b$, it holds
\begin{equation}
u_{2n} \sim C_{11}n^{-3/2}\rho^{-2n} \text{ as } n \to \infty
\label{1.18f}
\end{equation}
for some constant $C_{11} > 0$ (apply \cite{L}, \cite {Sa} to the random walk
$\{S_{2n}\}_{n \ge 0}$). Actually $C_{11}$ can be
explicitly given, but its precise value has no importance for us.
Clearly \eqref{1.18f} implies \eqref{1.18e} and \eqref{1.18ea}
and hence \eqref{1.18b} with $p=2$ (the period
for $\{S_n\}$).

We shall show after the proof of Theorem 3 that
\begin{equation}
F = \frac 1b,\; \rho = \frac{b+1}{2 \sqrt b} \text{ and } F(\rho)
= \frac{b+1}{2b}. \label{1.18j}
\end{equation}
Thus if $b \ge 2$, then in this example,
the limit in probability of $(1/2n)R_{2n}$
conditioned on $\cE_{2n}$ is strictly less than the almost sure limit
of $(1/n)R_n$ for the unconditioned walk. (These limits are
$(b-1)/(2b)$ and $(b-1)/b$, respectively.) Similar results hold for
random walks on various free products (see \cite{Ca}). For instance
we again have \eqref{1.18f}, and hence \eqref{1.18b}, for the simple
random walk on the free product of s copies of $\mathbb Z^q$ with $s
\ge 2, q \ge 1$, if we use the natural set of generators consisting
of the coordinate vectors in each factor $\mathbb Z^q$.

\bigskip
We end this section with some related remarks on the range and bridges
and list a number of open problems.

Note that if $\cG$ is a regular graph, that is all degrees are
equal, then the distribution of a simple random walk bridge is
just uniform measure on all $n$-step paths which return to the
starting point at the $n$-th step. The question of {\bf sampling}
a bridge on a given Cayley graph seems hard in general. We don't
even know how to sample a bridge on the lamplighter group. Or for
instance consider a simpler exercise, given a symmetric random
walk on $\mathbb Z$ which can take jumps of size 1 or 2. How can
one describe a bridge of such a random walk ?

\medskip \noindent
{\bf Open problem 1.} Find a necessary and sufficient condition for
\eqref{1.16} for a symmetric random walk on a Cayley graph $(G,\cS)$.
Theorems 1-3 suggest that perhaps amenability of $G$ is such a
necessary and sufficient condition (recall that under \eqref{1.5} and
\eqref{1.6} $G$ is amenable if and only if \eqref{1.10} holds; see
\cite{K1} or \cite{KV}, Section 5, or
\cite{Pa}, Theorems 4.19, 4.20 and Problem 4.24).

\medskip \noindent
{\bf Open problem 2.} Does there exist a transient random walk on a
Cayley graph or a simple random walk on a vertex transitive graph
for which $(1/n)R_n$ conditioned on $\cE_n$ tends to 0 in
probability ?
\newline
{\bf Open problem 3.} If the limit in probability of
$(1/n)R_n$, conditioned on $\cE_n$ exists, is it necessarily $\le
1-F$  (even if \eqref{1.10} fails) ?
\newline
\noindent {\bf Open problem 4.} Does unrestricted random walk
 drift further away from the starting point than a bridge ? More
 precisely, let $\cG$ be a vertex transitive graph, $\{S_k\}$ a simple
 random walk starting at a given vertex $e$, and
$\{S^b_k\}$, a simple random walk bridge conditioned to be back at
$e$ at time $n$. Is it true that for every fixed $k < n$, $d(S_k)
:=$ the (graph) distance of $S_k$ to $e$, stochastically dominates
$d(S^b_k)$ ? A related conjecture might be that on any vertex
transitive graph the simple random walk bridge is at most
diffusive, that is, $\max_{k \le n} d(S^b_k)$ is (stochastically)
at most of order $\sqrt n$ or perhaps $\sqrt n (\log n)^q$ for
some $q$. In a more quantative way one may ask for limit laws and
further bounds on the distribution of $\max_{k \le n} d(S^b_k)$
(see also Section 7 in \cite{E}). Note that on any graph with a
Gaussian off diagonal correction, i.e., for which there is a bound
$P\{\max_{k \le n} d(S^b_k) = \l\} \sim
C_{12}n^{-q}\exp[-C_{13}\l^2/n]$ for some positive $q$ and
$C_{12}, C_{13}$, the bridge is at least diffusive. Vertex
transitivity of $\cG$ seems to be crucial as it is not hard to
build examples of graphs on which bridges have non-diffusive and
more erratic behavior.

A weaker statement might be that on any vertex transitive graph
there is a subexponential lower bound of the form $P\{d(S^b_k) \le
\l|\cE_n\} \le \exp[-C\l/n]$ for the conditional probability given
$\cE_n$, that the bridge is near the starting point at time $k$.
\comment{
These conjectures or problems
can be considered in the context of general graphs, not just Cayley
graphs or vertex transitive graphs.
endcomment}


Another interesting problem in the same direction is to adapt the
Varopoulos-Carne subgaussian estimate (see \cite{LP}, Theorem 12.1)
to bridges. That is, to prove that $P(S^b_k =v) \leq e^{-d^2(v)/2n}$
or a weaker result of this type. In \cite{K2}, Proposition 3.3, there
is a proof that can give a subgaussian estimate for bridges on
certain graphs.

\medskip
We add two more observations which may be of interest, but do not
merit the label ``open problem." The first concerns a random walk on
the lamplighter group $\mathbb Z_2 \wr \mathbb Z^d$. A generic
element of the group is of the form $(\si, y)$ with $y
\in \mathbb Z^d$ and $\si$ a function from $\mathbb Z^d$ into $\{0,1\}$.
Assume that the random walk starts in $(\si_0,\bold 0)$, where
$\si_0$ is the zero function and $\bold 0$ is the origin in $\mathbb
Z^d$. Further, let $\si_y$ be the configuration obtained from $\si$
by changing $\si(y)$, the value of $\si$ at the position $y$, to
$\si(y) + 1 \mod 2$.

Assume that the random walk moves from $(\si,y)$ to $(\si, y \pm
e_i)$ or to $(\si_y, y\pm e_i),\; 1
\le i \le d$, with probability $1/(4d)$ for each of the possibilities. Here
$e_i$ is the $i$-th coordinate vector in $\mathbb Z^d$. Let us define
$\si_k,y_k,\si_k^b, y_k^b$ by $S_k = (\si_k,y_k)$ and $S_k^b
=(\si_k^b, y_k^b)$. Then, for the unrestricted random walk, $y_k$ is
simple random walk on $\mathbb Z^d$. For the bridge, all possible
paths have equal weight. Let $(y_0 = \bold 0, y_1, \dots, y_n =
\bold 0)$ be a nearest neighbor path on $\mathbb Z^d$ from $\bold 0$
to $\bold 0$. How many sequences $\{\si_k, y_k\}_{0\le k \le n}$ are
there which return to $(\si_0,\bold 0)$ at time $n$ and which project
onto the given sequence $\{y_k\}$ ? Such a sequence must have its
first coordinate set to 0 at position $y$ at the last visit to $y$ by
$\{y_k\}_{0 \le k \le n}$. Moreover, $\si_k$ can change only at the
position $y_k$. If $y$ is not one of the $y_k$ there is no condition
on the $\si_k$ at position $y$ at all. If $N_n$ denotes $|\{y_0, y_1,
\dots y_{n-1}\}|$ (i.e., the range of the projection on $\mathbb Z^d$
of the random walk on $\mathbb Z_2 \wr \mathbb Z^d$), then one sees from
the above that the number of possible $\{\si_k, y_k\}_{0\le k \le n}$
which return to $(\si_0,\bold 0)$ at time $n$ and which project onto
the given sequence $\{y_k\}$ equals $2^{n-N_n}$. Thus, if $p_r$
denotes the probability that an $n$-step simple random walk bridge
from $\bold 0$ to $\bold 0$ in $\mathbb Z^d$ has range $r$, then the
probability that the projection of the random walk bridge on $\mathbb
Z_2 \wr \mathbb Z^d$ has range $r$ equals
\begin{equation*}
P\{N_n = r\} = \frac {2^{n-r}p_r}{\sum_s 2^{n-s}p_s} = \frac {2^{-r}p_r}{\sum_s
2^{-s}p_s},
\end{equation*}
(this formula is very similar to equation (3.1) in \cite{PS})
and the expectation of the range of the projection is
\begin{equation}
E\{N_n\} = \frac {\sum_{1 \le r \le n} r2^{-r}p_r}
{\sum_{1 \le s \le n}2^{-s}p_s}.
\label{1.18k}
\end{equation}

We shall next argue that this expectation is only $O(n^{d/(d+2)})$,
so that by time $n$ the projection on $\mathbb Z^d$ of the bridge of
the random walk only travels a distance $O(n^{d/(d+2)})$ from the
origin. If $d=1$ it can be shown that then also the bridge of the
random walk on $\mathbb Z_2 \wr \mathbb Z$ itself only moves
distance $O(n^{1/3})$ from the starting point by time $n$. More
precisely, define for the bridge
\begin{equation*}
D_n = \max_{0 \le k \le n} d(S^b_k).
\end{equation*}
Then $n^{-1/3}D_n, n \ge 1$, is a tight family.
Even though this maximal distance $D_n$ for the bridge
grows slower than $n$, the range of the bridge
still grows in first order at the same rate as the range of the
unrestricted random walk (by example (iii) to Theorem 2).
(We note in passing that a similar, but weaker, comment applies if $\{S_n\}$
is simple random walk bridge on the regular $b$-ary tree $\cG_b$.
In this case it is known (\cite{BJ,NGV}) that $n^{-1/2}D_n, n \ge 1$,
is a tight family. Now the range of the bridge is in first order
still linear in $n$, but the ranges of the bridge and of
the unrestricted  simple random walk on $\cG_b$ differ already in
first order; see the example to Theorem 3). We point out that
\cite{E} argues that $n^{-1/3}E\{D_n\} \ge n^{-1/3}d(S_{n/2})$ is also
bounded away from 0 for $n$ even, $d=1$ and a random walk which differs
from ours only a little in the choice of the distribution of the
$X_i$. Presumably such a lower bound for $D_n$ also holds in our case,
but is not needed for the present remark.

Here is our promised estimate for $E\{N_n\}$.
It is known (see \cite {DoV}) that for an unrestricted simple random
walk $\{S_k\}$ on  $\mathbb Z^d$,
\[
\begin{split}
P\{R_n \le (2L+1)^d\} &\ge P\{S_k \in [-L, L]^d, 0 \le k \le n\}\\
& \ge
C_{14}\exp[-C_{15}nL^{-2}]
\end{split}
\]
for some constants $0 < C_i < \infty$.
It can be shown from this that also
\[
P\{S_k \in [-L, L]^d, 0 \le k \le n, S_n = \bold 0\}
\ge C_{16}L^{-d}\exp[-C_{15}nL^{-2}].
\]
A fortiori
\[
P\{S_k \in [-L, L]^d, 0 \le k \le n\big| S_n = \bold 0\} \ge
C_{16}L^{-d}\exp[-C_{15}nL^{-2}].
\]
By taking $L = n^{1/(d+2)}$ we find that the denominator in
\eqref{1.18k} is for large $n$ at least
\[
2^{-(2L+1)^d }C_{16}L^{-d}\exp[-C_{15}nL^{-2}]
\ge C_{17}\exp[-C_{18}n^{d/(d+2)}].
\]
Consequently the right hand side of \eqref{1.18k} is at most
\[
\begin{split}
&\frac{\sum_{r \le K} r2^{-r}p_r}{\sum_s 2^{-s}p_s} +
\frac n{C_{17}\exp[-C_{18}n^{d/(d+2)}]}
  \sum_{K < r \le n}2^{-r}p_r\\
&\le K + \frac n{C_{17}}\exp[C_{18}n^{d/(d+2)}] 2^{-K}.
\end{split}
\]
By choosing $K = 2C_{18}n^{d/(d+2)}/ \log 2$ we obtain the promised
bound $E\{N_n\} = O(n^{d/(d+2)})$.

\smallskip
\bigskip

Our last observation deals with bridges on special finite graphs.
Let $\cG$ be a $\cD$-regular vertex transitive {\bf expander} of
size N with girth (smallest cycle) of size $c \log N$. The mixing
time for simple random walk on such a graph is $C \log N$, for
some constant $C > c$, where mixing time is taken in the strong
sense of the maximum relative deviation. That is, we take the
mixing time to be the number of steps, $k$, a simple random walk
has to take to make $\sup_{v,w \in \cG}\big|P\{S_k = w|S_0 =
v\}/\pi(w) -1\big|$ smaller than some prescribed number, where
$\pi(\cdot)$ is the stationary measure for the random walk (see
\cite{Si}). Now a bridge of length $n < c\log N$ is just a bridge
on a $(\cD-1)$-regular tree, since there are no cycles of length
$< c\log N $. Thus the range of such a bridge is in first order
described by Theorem 3. On the other hand, a bridge of length $n>
C'\log N$ for large enough $C'$ may be expected to look like an
unconstrained simple random walk at least for times in $[c\log N,
n-c\log N]$ , because of the short mixing time. For larger $n$,
but still order $\log N$, we actually expect the range of such a
bridge to be like the range of a $n$-step unconstrained random
walk on a $(\cD-1)$-regular tree, that is $n(\cD-2)/(\cD-1)$ in
first order. This suggests that maybe there is a critical $C^*$,
so that a bridge of length $< C^* \log N$ (respecively $> C^*\log
N$) looks as in the first case (respectively, second case).

\medskip
\noindent {\bf Acknowledgement.} We thank Laurent Saloff-Coste for
several helpful conversations about the subject of this paper.

\section{{\bf Proofs.}}
{\it Proof of Theorem 1.} We shall treat the aperiodic case (i.e.,
the case with $p =1$) only. The case when the period equals 2 can
be treated in the same way. One merely has to restrict all the
subscripts to even integers.

Define the random variables
$$
Y(k,M) : = I\big[S_{k+r} \ne S_k \text{ for } 1 \le r \le M\big].
$$
By a last exit decomposition we have for any positive integer $M$
\begin{equation}
\begin{split}
R_{n} &= \sum_{k=0}^{n-1} I\big[S_{k+r} \ne S_k \text{ for } 1 \le
  r \le n-k\big]\\
&\le M + \sum_{k=0}^{n-M} I\big[S_{k+r} \ne S_k \text{ for } 1 \le
  r \le M\big]\\
&= M + \sum_{k=0}^{n-M} Y(k,M)
\end{split}
\label{2.5}
\end{equation}
(compare proof of Theorem 4.1 in \cite{Sp1}, which uses a
first entry decomposition). Now by the Markov property of the
random walk $\{S_k\}$ and the transitivity of $\cG$,
\begin{equation}
\begin{split}
&P\{Y(k,M) =1\big |S_0, S_1,\dots, S_k\}\\
&= P\{S_{k+r} \ne S_k \text{ for } 1 \le r \le M \big|S_0, S_1,
\dots,
S_k\}\\
& = P\{S_r \ne S_0 \text{ for } 1 \le r \le M\} = P\{Y(0,M) =1\}.
\end{split}
\label{1.20}
\end{equation}
This relation says that any collection $\{Y(k_i,M)\}$ of these
random variables with $|k_i - k_j| > M$ for $i \ne j$ consists of
i.i.d. random variables. In fact, each $Y(k,M)$ can take only the
values 0 or 1. This will allow us to use exponential bounds for
binomial random variables.

First observe that for given $\ep > 0$ we can choose $M$ such that
\begin{equation}
\begin{split}
E\{R_n\} &\ge \sum_{k=0}^{n-1} P\{S_{k+r} \ne S_k \text{ for all
}r
\ge 1\}\\
& \ge \sum_{k=0}^{n-1} \big[E\{Y(k,M)\} - P \{\text{first
    return to $S_k$ occurs at time}\\
&\phantom{MMMMMMMMMMMMMM} k+r \text{ for some } r > M\} \big]\\
&=\sum_{k=0}^{n-1} \big[E\{Y(k,M)\} - \sum_{r > M} f_r\big] \ge
\sum_{k=0}^{n-1} E\{Y(k,M)\} - \ep n.
\end{split}
\label{2.5a}
\end{equation}
We now rewrite \eqref{2.5} as
\begin{equation}
R_n \le M+\sum_{a=0}^M \sum_{\substack{k \equiv a \text{ mod }
(M+1)\\0 \le k \le n-1}} Y(k,M).
\end{equation}
Moreover, by \eqref{2.5} we have for large enough $M$ and $n \ge
n_0$ for some $n_0 = n_0(M)$,
\begin{equation}
\begin{split}
E\{R_n\} & \le M + \sum_{k=0}^{n-M} E\{Y(k,M)\} \\
&\le M + nP\{S_r \ne S_0 \text{ for } 1\le r\le M\} \le n[1-F
+\ep].
\end{split}
\label{2.5b}
\end{equation}
Thus, if for given $\ep > 0$, $M$ is chosen so that \eqref{2.5a}
and \eqref{2.5b} hold, then for all large $n$
\begin{equation}
\begin{split}
&P\{R_n \ge n(1-F) + 3\ep n\} \le
P\big\{R_n - E\{R_n\} \ge 2\ep n\big\}\\
&\le \sum_{a=0}^M P\Big\{  \sum_{\substack{k \equiv a \text{ mod }
(M+1)\\0 \le k \le n-1}} \big[Y(k,M)- E\{Y(k,M)\}\big] \ge \frac
\ep{M+1} n\Big\}.
\end{split}
\end{equation}
The right hand side here tends to 0 exponentially fast as $n \to
\infty$ by standard exponential bounds for large deviations in
binomial distributions (e.g. by Bernstein's inequality
\cite{ChT}, Exercise 4.3.14).

We now remind the reader of \eqref{1.10}. It is easy to
see from the definition of $u_k$ that
\begin{equation}
u_{2k+2\l} \ge u_{2k}u_{2\l}. \label{2.5d}
\end{equation}
 From this and the fact that
$u_{2k} \ge [u_2]^k >  0$ one obtains that $\limn [u_{2n}]^{1/n}$
exists, and then by \eqref{1.10} that this limit must equal 1. This
holds regardless of whether $p=1$ or 2. If $p=1$ then even $u_m >
0$ for some odd $m$ and $u_n \ge u_{n-m} u_m$ then shows that also
$[u_{2k+1}]^{1/(2k+1)} \to 1$ as $k \to \infty$. Consequently,
$$
\frac 1{u_n} P\{\frac 1n R_n > 1-F+3\ep\} \to 0 \text{ as } n \to
\infty.
$$
(in fact, exponentially fast). But then also
$$
P\{\frac 1n R_n > 1-F+3\ep\big|\cE_{2n}\} \to 0.
$$
This proves \eqref{1.11}. If the random walk $\{S_n\}$ is
recurrent, then $F =1$ en \eqref{1.11} says that $(1/n)R_n$
conditioned on $\cE_n$ tends to 0 in probability. \qed

\bigskip \noindent
{\it Proof of Theorem 2.} We begin this subsection with
 some lemmas on the smoothness of the $u_n$ (as functions of $n$).
Then we give a general sufficient condition in terms of the $u_n$
and $f_n$ for the limit in probability of $(1/n)R_n$ conditioned
on $\cE_n$ to equal $1-F$. Finally we show that this sufficient
condition holds under the conditions of Theorem 2.
\begin{lm} Let $\{S_n\}$ be a random walk on  a Cayley graph $(G,\cS)$
  and assume that the symmetry property \eqref{1.5} holds. Then there
  exists a probability measure $\mu$ on $[-1, 1]$ such that
\begin{equation}
u_n = \int_{[-1,1]} x^n \mu(dx), \;n \ge 0. \label{2.45}
\end{equation}
Consequently
\begin{equation}
u_{2n} \text{ is non-increasing in $n$}, \label{2.46}
\end{equation}
\begin{equation}
u_{2n+1} \le u_{2n}, \label{2.46a}
\end{equation}
and also
\begin{equation}
u_{2r}u_{2n-2r}\text{ is non-increasing in $r$ for } 0 \le r \le
\frac {n-1}2 \le \frac {n+1}2 \le n-r. \label{5.38}
\end{equation}
\end{lm}

\noindent{\it Proof.} We write
$$
P(v,w) = P\{X = v^{-1}w\} = P\{S_{n+1} = w|S_n= v\}, \quad v,w \in
G,
$$
for the transition probabilities of the random walk $\{S_n\}$.
Then the $k$-th power of $P$ gives the $k$-step transition
probability. That is
$$
P^k(v,w) = P\{S_{k+n} = w|S_n=v\}.
$$
$P$ defines a linear operator on $\l^2(G)$ by means of
$$
Pf(v) = \sum_{w \in G}P(v,w)f(w).
$$
This linear operator takes $\l^2(G)$ into itself and is
self-adjoint. By the spectral theorem (see \cite{R}, Theorems
12.23, 12.24; see also \cite{KV}, Section 5)
there therefore exists a measure $\mu$ on the Borel
sets of $\Bbb R$ such that for $f_0(v) = I[v = e]$, and
$\langle \cdot, \cdot \rangle$ the inner product on $\l^2(\cG)$,
$$
u_n  = \langle f_0, P^nf_0\rangle = \int_{\Bbb R} x^n \mu(dx), \; n \ge 0.
$$
From $|Pf(v)|^2\le \sum_w P(v,w)|f|^2(w) \sum_w P(v,w) = \sum_w
P(v,w)|f|^2(w)$ we see that $\|P\| \le 1$, so that the support of
$\mu$ must be contained in $[-1,1]$ and \eqref{2.45} must hold.
(This can also be see directly from $|u_n| \le 1$ for all $n$.)

This proves the existence of some measure $\mu$ for which
\eqref{2.45} is satisfied. In fact, the spectral theorem tells us
that for $A$ a Borel set of $[-1,1],\; \mu(A) =
\langle f_0,E(A)f_0\rangle$ for some resolution of the identity $\{E(\cdot)\}$.
In particular, $E([-1,1])$ is the identity $I$ on $\l^2(\cG)$ so that
$\mu([-1,1]) = \langle f_0,If_0\rangle = 1$. Thus $\mu$ is a
probability measure as claimed.

\eqref{2.46} is an immediate
consequence of \eqref{2.45}. As for \eqref{5.38}, we have from
\eqref{2.45} that
\begin{equation}
\begin{split}
&u_{2r}u_{2n-2r} - u_{2r+2}u_{2n-2r-2}\\
&=\int \mu(dx)\int\mu(dy) \big[x^{2r}y^{2n-2r} -
x^{2r+2}y^{2n-2r-2}\big]\\
&= \frac 12 \int \mu(dx)\int\mu(dy) \big[x^{2r}y^{2n-2r} -
x^{2r+2}y^{2n-2r-2}\\
&\phantom{MMMMMMMMMMMM}+y^{2r}x^{2n-2r}
-y^{2r+2}x^{2n-2r-2}\big]\\
&=\frac 12 \int \mu(dx)\int\mu(dy) x^{2r}y^{2r}\big[y^{2n-4r-2}
-x^{2n-4r-2}\big] \big[y^2 -x^2\big] \ge 0,
\end{split}
\end{equation}
because $[y^{2n-4r-2} -x^{2n-4r-2}][y^2 -x^2] \ge 0$ for all $x,y$
if $2n-4r-2 \ge 0$. \qed

\begin{lm}
Let $\{S_n\}$ be a random walk on  a Cayley graph $(G,\cS)$ which
satisfies the symmetry assumption \eqref{1.5}.
Assume that
\begin{equation}
u_{2n} \ge e^{-g(n)} \label{2.25}
\end{equation}
for some function $g(\cdot)$ which satisfies
\eqref{1.12}, \eqref{1.13}.
Then
\begin{equation}
\frac n{n+r} \Big(\frac{r}{n+r}\Big)^{r/n} e^{[-rg(n)/n]}
 \le \frac {u_{2n+2r}}{u_{2n}}
\le 1.
\label{2.23}
\end{equation}
for all $r \ge 1$. In particular,
\begin{equation}
\limn \frac{u_{2n+2}}{u_{2n}} = 1.
\label{2.23aa}
\end{equation}
Moreover, if $p=1$,
\comment{
there exists a further odd integer $L_0 \ge 1$
such that
\begin{equation}
\begin{split}
&\big(1+o_n(1)\big)
 \exp\Big[-C_7 \sum_{j=0}^{r+(L_0-1)/2}\sqrt{\frac
{g(\lceil n/2\rceil+j)}{n+2j}}\Big]\\
&\le\frac{u_{n+2r+1}}{u_n} = \frac {u_{n+2r+1}} {u_{n+2r+L_0+1}}
\cdot
\frac {u_{n+2r+L_0+1}}{u_n} \\
&\le  \big(1+o_n(1)\big)
 \exp\Big[C_7 \sum_{j=0}^{r+(L_0-1)/2}
\sqrt{\frac{g(\lceil n/2\rceil+j)}{n+2j}}\Big],
\end{split}
\label{2.37}
\end{equation}
where $o_n(1) \to 0$ is a quantity which tends to 0 as $n \to
\infty$. In particular
\endcomment}
\begin{equation}
\limn \frac {u_{n+p}}{u_n} = 1. \label{2.40}
\end{equation}
\end{lm}
\noindent {\it Proof.}
The right hand inequality in \eqref{2.23} is part of Lemma 1.

To find a lower bound for $u_{2n+2r}/u_{2n}$
we again appeal to \eqref{2.45}. This tells us that for any $\ga \ge 0$,
\begin{equation*}
\begin{split}
e^{-g(n)}&\le
u_{2n} =\int_{[-1,1]} x^{2n}\mu(dx)\\
&=  \int_{|x| \le \exp[-(\ga +
      g(n))/(2n)]} x^{2n}\mu(dx) + \int_{|x| > \exp[-(\ga +
      g(n))/(2n)]} x^{2n}\mu(dx)\\
&\le \exp[-(\ga+g(n)) + \int_{|x| > \exp[-(\ga +
      g(n))/(2n)]} x^{2n}\mu(dx)\\
&\le e^{-\ga}u_{2n} +  \int_{|x| > \exp[-(\ga +
      g(n))/(2n)]} x^{2n}\mu(dx).
\end{split}
\end{equation*}
Consequently,
\begin{equation*}
[1-e^{-\ga}] u_{2n} \le \int_{|x| > \exp[-(\ga +
      g(n))/(2n)]} x^{2n}\mu(dx).
\end{equation*}
It follows that
\[
\begin{split}
u_{2n+2r} &\ge \int_{|x| > \exp[-(\ga + g(n))/(2n)]} x^{2n+2r}\mu(dx)\\
&\ge \exp\big[-\frac{r}n (\ga +g(n))\big]
\int_{|x| > \exp[-(\ga + g(n))/(2n)]} x^{2n}\mu(dx)\\
&\ge \exp\big[-\frac{r}n (\ga +g(n))\big] [1-e^{-\ga}]u_{2n}.
\end{split}
\]
The left hand inequality of \eqref{2.23} follows by taking $e^{-\ga} =
r/(n+r)$.

The limit relation \eqref{2.40} is proven in \cite {Ge} and in \cite{Av}.
\qed

\begin{lm}
Assume that $\{S_n\}$ is a random walk on an infinite Cayley graph
for which \eqref{1.5}, \eqref{1.6} and \eqref{1.10} hold. If, in
addition, for each $\eta > 0$
\begin{equation}
\lim_{M \to \infty} \limsup_{n \to \infty,p|n}\;
\sum_{\substack{M\le r \le (1-\eta)n\\p|r}}
 f_r\frac {u_{n-r}}{u_n} = 0,
\label{5.52}
\end{equation}
then for all $\ep > 0$
\begin{equation}
\lim_{n \to \infty, p|n} P\{\big|\frac 1n R_n -1+F\big| > \ep \big|\cE_n\}
=0. \label{5.53}
\end{equation}
\end{lm}
\noindent{\it Proof.} We claim that it suffices to show
\begin{equation}
\liminf_{n \to \infty, p|n} E\{\frac 1n R_n\big|\cE_n\} \ge 1-F.
\label{2.21}
\end{equation} Indeed, we already know from Theorem 1 that
\eqref{1.11} holds. Together with $(1/n)R_n \le 1$ this implies
for $\ep > 0$ and $n$ a large multiple of $p$ that
\begin{equation}
\begin{split}
1-F - \ep &\le E\{\frac 1n R_n\big|\cE_n \} \le \big(1-F- \sqrt
\ep\big)
P\big\{\frac 1nR_n \le 1-F-\sqrt \ep\big|\cE_n\big\}\\
&\phantom{MM} + \big(1-F+ \ep\big)
P\big\{1-F - \sqrt \ep < \frac 1n R_n \le  1-F+\ep\big|\cE_n \big\}\\
&\phantom{MM}+P\big\{\frac 1n R_n > 1-F+\ep\big|\cE_n \big\}\\
&\le \big(1-F- \sqrt \ep\big)
P\big\{\frac 1nR_n \le 1-F-\sqrt \ep\big|\cE_n\big\}\\
&\phantom{MM} + \big(1-F+ \ep\big)\big[1- P\big\{\frac 1nR_n \le
1-F-\sqrt \ep\big|\cE_n\big\}\big] + \frac 12 \ep.
\end{split}
\label{2.22}
\end{equation}
By simple algebra this is equivalent to
$$
P\big\{\frac 1n R_n \le 1-F - \sqrt \ep\big|\cE_n\big\} \le \frac
{5\ep}{2(\sqrt \ep + \ep)}.
$$
This justifies our claim that we only have to prove \eqref{2.21}.

We turn to the proof of \eqref{2.21}. Again by a last exit decomposition
\begin{equation*}
\begin{split}
R_n &\ge R_{\lfloor(1-\eta)n\rfloor}= \sum_{0 \le k <(1-\eta)n}
I[S_{k+r} \ne S_k \text{ for } 1 \le r \le (1-\eta)n-k]\\
&=\sum_{0 \le k <(1-\eta)n} \big[1 - \sum_{r=1}^{(1-\eta)n -k}
I[\text{first return to $S_k$ occurs at time $k+r$}]\big]
\end{split}
\end{equation*}
Now multiply this inequality by $I_n := I[\cE_n]$ and take
expectations. This yields
\begin{equation}
\begin{split}
\frac 1n E\{&R_n|\cE_n\} = \frac 1{nu_n} E\{R_nI_n\} \ge \frac{(1-\eta)n}n \\
&-\frac 1{nu_n} \sum_{0\le k <(1-\eta)n} \;\sum_{1 \le r <
(1-\eta)n-k}\;
\sum_{v \in G} P\{S_k = v, \text{first return to $v$}\\
&\phantom{MMMMMMM}\text{after $k$ is at time $k+r$ and }S_n = e\}.
\end{split}
\label{5.51}
\end{equation}
The inner sum of the triple sum in the right hand side here equals
\begin{equation}
\begin{split}
\sum_{v \in G}& P\{S_k = v\} f_r P\{S_n=e|S_{k+r} = v\} \\
&=f_r\sum_{v \in G} P\{S_k=v\}P\{S_{n-k-r} = e|S_0 = v\} =
f_ru_{n-r}.
\end{split}
\label{5.51a}
\end{equation}
We substitute this into \eqref{5.51} and use the assumption
\eqref{5.52} to obtain for any $\eta \in (0,1)$
\begin{equation}
\liminf_{n \to \infty, p|n} \frac 1n E\{R_n|\cE_n\} \ge 1-\eta -
\limsup_{M
  \to  \infty} \limsup_{n \to \infty, p|n}  \sum_{r=1}^M f_r
  \frac{u_{n-r}}{u_n}.
\label{5.55}
\end{equation}
Note that $f_r\le u_r = 0$ if $p \nmid r$. Thus, if we prove that
\begin{equation}
\lim_{n \to \infty, p|n} \frac {u_{n-r}}{u_n} = 1 \text{ for fixed
$r$
  with $p|r$},
\label{5.56}
\end{equation}
then the desired \eqref{2.21}, and hence also \eqref{5.53} will
follow. We shall now deduce \eqref{5.56} from Lemma 2. As we saw
in \eqref{2.5d} and the lines following it, if $p=1$, then
\eqref{1.10} implies that
\begin{equation}
\limn [u_{2n}]^{1/n} = 1. \label{5.50}
\end{equation}
Define
\[
g(n) = \log \frac 1{u_{2k}}.
\]
$g(\cdot)$ is nondecreasing by virtue of \eqref{2.46},
and its limit as $n \to \infty$ must be $\infty$. To
see this note that the Markov chain $\{S_n\}$ cannot be positive
recurrent, because the measure which puts mass 1 at each vertex of
$\cG$ is an infinite invariant measure if $G$ is infinite (see
\cite{F}, Theorem in Section XV.7 and Theorem XV.11.1). Thus it
must be the case that $u_n \to 0$. It follows that $g(\cdot)$
satisfies \eqref{1.12}. Moreover, \eqref{1.13} for $g(\cdot)$ is
implied by \eqref{5.50}, and \eqref{2.25} holds
by definition. \eqref{5.56} then follows from \eqref{2.23aa} and
\eqref{2.40}.
\qed

We are finally ready to give the proof of Theorem 2. By virtue of
Lemma 3 it suffices to prove that \eqref{1.10} and \eqref{5.52}
hold. Assume first that \eqref{1.16a} holds. \eqref{1.10} then
follows easily. Indeed, \eqref{1.16a} implies that there exists
some $n_0$ and a constant $C_{14} > 0$ such that
\begin{equation}
u_{2^rn_0} \ge C_{14}^{r-1}u_{2n_0}.
\label{2.93}
\end{equation}
From the existence of $\limn [u_{2n}]^{1/n}$ (see \eqref{2.5d} and
the
  lines following it) we then see that \eqref{1.10} holds. As for
  \eqref{5.52}, we have from \eqref{1.16a}
and the monotonicity in $q$ of $u_{2q}$ that for fixed $\eta > 0$
there exists some constant $C_{15}$ depending on $\eta$ only such that
$u_{n-r}/u_n \le C_{15}$ for even $n$ and even $r \le (1-\eta)n$. For
the case when $p=2$ this shows for even $n$ that
\[
\sum_{\substack{M \le r< (1-\eta)n\\ r \text{ even}}} f_r\frac
 { u_{n-r}}{u_n} \le C_{15}\sum_{r \ge M}f_r.
\]
In this case \eqref{5.52} is therefore immediate. If $p=1$ we can
use essentially the same argument, since for odd $r \le
(1-\eta)n$, $u_r \le u_{r-1}$ (see \eqref{2.46a}), while for odd
$n$, $u_n \ge u_{L_0}u_{n-L_0}$, where $L_0$ is a fixed odd
integer for which $u_{L_0} > 0$. These simple observations show
that in case $p=1$ we still have $u_{n-r}/u_n$ bounded by some
$C_{15}(\eta)$ for all $r \le (1-\eta)n$, and hence also \eqref{5.52}
holds.

Now we turn to the proof of \eqref{1.16} from
\eqref{1.12}-\eqref{1.15}. First note that \eqref{1.13} and
\eqref{2.25} imply that \eqref{1.10} (or even \eqref{5.50}) holds
(see a few lines after \eqref{2.5d}). Thus it is again enough to
prove \eqref{5.52}.
To this end, define
\[
\tau = 2 \Big \lfloor \frac n{g(n/2)}\Big \rfloor
\label{5.62}
\]
and use Lemma 2 to obtain for fixed $M$, but large even $n$
\begin{equation}
\begin{split}
&\sum_{M\le r < (1-\eta)n}  f_r\frac{u_{n-r}}{u_n} \le \sum_{M \le
r \le n-M}
\frac {u_ru_{n-r}}{u_n}\\
& = 2\sum_{M\le r\le (n-1)/2} \frac {u_ru_{n-r}}{u_n}\\
&\le 2\sum_{M\le r\le (n-1)/2,2|r} \frac {u_ru_{n-r}}{u_n}
+2\sum_{M-1\le r\le (n-2)/2,2|r} \frac {u_ru_{n-2-r}}{u_n}
\text{ (by \eqref{2.46a})}\\
&\le 2\sum_{r=M}^{\tau-1} \frac {u_ru_{n-r}}{u_n} +2n \frac
{u_{\tau} u_{n-\tau}}{u_n}
+2\sum_{r=M-1}^{\tau-1}\frac{u_ru_{n-2-r}}{u_n}
+2n \frac {u_{\tau} u_{n-2-\tau}}{u_n} \text{ (by \eqref{5.38})}\\
&\le 4\sum_{r=M-1}^\infty u_r \frac{u_{n-2-\tau}}{u_n}+ 4n \frac
{u_{\tau} u_{n-2-\tau}}{u_n} \text{ (by \eqref{2.46} and
\eqref{2.46a})}.
\end{split}
\label{5.60}
\end{equation}
Next, \eqref{2.23} (with $2n$ replaced by $n-2-\tau$ and $2r$ by
$\tau+2$) shows that for large even $n$,
\begin{equation}
\begin{split}
&\frac {u_{n-2-\tau}}{u_n} \\
&\le
\frac n{(n-\tau-2)} \Big(\frac n{\tau +2}\Big)^{(\tau+2)/(n-\tau + 2)}
\exp\big[\frac {\tau+2}{n-\tau-2}g\big((n-\tau-2)/2\big)\big].
\end{split}
\label{5.60ab}
\end{equation}
But the definition of $\tau$ and the monotonicity of $g(\cdot)$ show
that $\tau = o(n)$ and
\[
\frac{\tau+2}{n-\tau -2}g\big((n-\tau-2)/2\big) \le 2\frac\tau n g(n/2)
\le 4
\]
for large $n$. Therefore, the right hand side of \eqref{5.60ab} is
bounded by some constant $C_{16}$.
Substitution of this bound into \eqref{5.60} now shows that
\begin{equation}
\sum_{M \le r < (1-\eta)n} f_r \frac{u_{n-r}}{u_n} \le 4C_{16}
\sum_{r=M-1}^\infty u_r + 4C_{16}nu_\tau. \label{5.61}
\end{equation}
Theorem 1 already proves \eqref{1.16} if $\{S_n\}$ is recurrent,
so we may assume that this random walk is transient. In this case
we can make the first term in the right hand side of \eqref{5.61}
small by choosing $M$ large (see \cite {F}, Theorem XIII.4.2).
Finally, by \eqref{1.15}, \eqref{5.62} and \eqref{1.14} also
\[
nu_\tau \le nh(\tau/2) = n h\big(\big\lfloor \frac n{g(n/2)}
\big\rfloor\big)\to 0.
\]
This completes the proof of \eqref{5.52} for even $n$. The case of
odd $n$, which arises only when $p = 1$ can be reduced to the
case of even $n$, because $u_{n-r}/u_n \sim u_{n+1-r}/u_{n+1}$ if
$p=1$, by virtue of \eqref{2.40}.
\qed

\bigskip
\noindent
{\it Proof of Theorem 3.} The properties in \eqref{1.18c} are well
known (see \cite{V}, Theorem C and also Theorem 4.1 in \cite{L}). We merely
add a few comments concerning \eqref{1.18c} which will be needed
for the examples.

As usual, if $\{S_n\}$ is a random walk on a Cayley graph
$(G,\cS)$ then we take $e$ to be the identity element of $G$. If
$\{S_n\}$ is simple random walk on some other vertex transitive
graph $\cG$ then $e$ is any fixed vertex of $\cG$. $f_n$ and $u_n$
are as in \eqref{1.5a}. In all cases $u_{n+m} \ge u_nu_m$. As we
already pointed out in \eqref{2.5d} and the following lines this
implies that $\limn [u_{np}^{1/(np)}]$ exists and is at least
equal to $u_2 >
  0$. Moreover, $u_m = 0$ if $p \nmid m$. Thus $U(z) = \sum_{n\ge
    0}u_{np}$ has the radius of convergence
$\rho = \big[\limn [u_{np}^{1/(np)}\big]^{-1} \in[1, \infty)$
(recall
      $|u_n| \le 1$). From the theory of recurrent events (see
      \cite{F}, Section XIII.3)
it then follows that
\[
U(z) = \frac 1{[1-F(z)]} \text{ for } |z| < 1,
\]
where $F(z) = \sum_{n=1}^\infty f_nz^n$, as defined
in the statement of the theorem. But the right hand side here is
analytic on the open disc $\{z: |z| <\rho\}$ with $\rho =
\min\{\rho_1, \rho_2\}$, where
$\rho_1$ equals the first
singularity of $F(\cdot)$ on the positive real axis, and $\rho_2 =
\sup \{x > 0: F(x) < 1\}$. Thus, also the power series $U(z) =
\sum_{n=0}^\infty u_nz^n$ converges at least for $|z| < \rho$. On the
other hand it cannot be that the powerseries for $U$ converges on all
of the disc $\{|z| < \rho + \ep\}$ for some $\ep >0$, because
$1/[1-F(z)]$ cannot be analytic in such a disc. In fact, since $F$ is
a powerseries with non-negative coefficients its smallest singularity
must be on the positive real axis. Thus the radius of convergence of
$U$ equals $\rho$. The fact that $\rho \le \rho_2$ (and Fatou's lemma)
shows that $F(\rho)\le 1$, as claimed in \eqref{1.18c}.

To prove \eqref{1.18b} from the conditions \eqref{1.18e} and
\eqref{1.18ea} we shall show that for all $0 < \eta < 1$
\begin{equation}
\lim_{n \to \infty,p|n} \frac 1n E\{R_{(1-\eta)n}\big|\cE_n\} =
(1-\eta)[1-F(\rho)],
\label{1.19a}
\end{equation}
and then show that a good approximation to $(1/n)R_{(1-\eta)n}$ on
$\cE_n$ has a conditional variance, given $\cE_n$, which tends to 0
as $n \to \infty$ (see \eqref{1.19h}). Since the approximation can
be made as precise as desired, and since $|R_n
-R_{(1-\eta)n}| \le \eta n$ this will give us \eqref{1.18b}.
(Here and in the sequel we drop the largest integer symbol in
$\lfloor (1-\eta)n \rfloor$ for brevity).

The proof of \eqref{1.19a} is straightforward. We define
 \begin{equation}
\begin{split}
&J(v, k, r)\\
&\qquad = I\big[S_k = v, \text{ first time after $k$ at which
$S_.$
    returns to $v$ is $k+r$}\big].
\end{split}
\end{equation}
Again, by a last exit decomposition we have
\begin{equation}
R_{(1-\eta)n} =  \sum_{0 \le k < (1-\eta)n} \sum_{v \in \cG}
\big[1- \sum_{1 \le r
  \le (1-\eta)n - k} J(v,k,r)\big]
\end{equation}
For brevity we write $W$ for the triple sum
\[
 \sum_{0 \le k < (1-\eta)n} \sum_{v \in \cG}
\sum_{1 \le r  \le (1-\eta)n - k} J(v,k,r).
\]
With $I_n = I[\cE_n]$ as before, we have as in \eqref{5.51a}
\begin{equation*}
E\{WI_n\} = \sum_{0 \le k < (1-\eta)n}\; \sum_{1
\le r  \le (1-\eta)n - k} f_ru_{n-r},
\end{equation*}
and hence
\begin{equation}
\begin{split}
&E\frac 1n \big\{R_{(1-\eta)n}\big|\cE_n\big\} = (1-\eta) - \frac
1n \sum_{0 \le k < (1-\eta)n} \sum_{\substack{1 \le r
<(1-\eta)n-k\\p|r}}
f_r\frac{u_{n-r}}{u_n}\\
& = (1-\eta)- \frac 1n \sum_{0 \le k < (1-\eta)n} \sum_{\substack{1 \le r
<(1-\eta)n-k\\p|r}} f_r\rho^r \frac
{u_{n-r}\rho^{n-r}}{u_n\rho^n}.
\end{split}
\label{1.19c}
\end{equation}
Now, by \eqref{1.18e},\eqref{1.18ea} and the fact that $F(\rho) <
\infty$ we have,for any fixed $\eta' > 0$,
uniformly in $k \in [0,(1-\eta -\eta')n)$ that
\begin{equation}
\lim_{n \to \infty, p|n} \sum_{\substack{1 \le r
<(1-\eta)n-k\\p|r}} f_r\rho^r \frac {u_{n-r}\rho^{n-r}}{u_n\rho^n}
= \sum_{r \ge 1, p|r} f_r\rho^r = F(\rho). \label{1.19d}
\end{equation}
Moreover, for each fixed $\eta > 0$, and $k \le (1-\eta)n$,
\[
\sum_{\substack{1 \le r <(1-\eta)n-k\\p|r}} f_r\rho^r \frac
{u_{n-r}\rho^{n-r}}{u_n\rho^n}
\]
is bounded in $n$. \eqref{1.19a} is immediate from these
observations and \eqref{1.19c}.

We also obtain that for any given $\de \in (0,1)$ we can choose $M
= M(\de)$, independent of $n$, such that
\[
\frac 1n E\big\{\sum_{0 \le k < (1-\eta)n} \sum_{v \in \cG}
\sum_{M < r \le (1-\eta)n - k} J(v,k,r)\big|\cE_n\big\} \le \de.
\]
This leads us to write
\[
W_M =
 \sum_{0 \le k < (1-\eta)n} \sum_{v \in \cG}
\sum_{1 \le r  \le M \land \big((1-\eta)n - k\big)} J(v,k,r).
\]
With this notation we find
\begin{equation}
\lim_{n \to \infty, p|n} \frac 1nE\{W_M|\cE_n\} =
(1-\eta)\sum_{r=1}^Mf_r\rho^r \label{1.19f}
\end{equation}
and
\begin{equation}
\frac 1n E\big\{\big|R_{(1-\eta)n}- (1-\eta) + W_MI_n\big|
\big|\cE_n\big\} \le \de. \label{1.19g}
\end{equation}
Next we estimate $E\{W^2_M|\cE_n\}$ for a fixed $M$. From the
definition of $W_M$ itt follows that
\begin{equation}
\begin{split}
&E\{W_M^2I_n\}\\
&=
 \sum_{0 \le k < (1-\eta)n} \sum_{v \in \cG}
\sum_{1 \le r  \le M \land \big((1-\eta)n - k\big)}
 \sum_{0 \le \l < (1-\eta)n} \sum_{w \in \cG}
\sum_{1 \le s  \le M \land \big((1-\eta)n - k\big)}\\
&\phantom{MMMMMMMMMMMMMMMMMM}E\{J(v,k,r)J(w,\l,s)I_n\}
\end{split}
\label{1.19e}
\end{equation}
Now note that for fixed $\l, \sum_{w \in \cG} \sum_{1 \le s \le m}
J(w,\l,s)$ is a sum of indicator functions of disjoint  events and
is therefore bounded by 1. Thus the terms in the multiple sum in
\eqref{1.19e} with $|\l-k| \le M$ contribute at most
\[
\begin{split}
& \sum_{0 \le k < (1-\eta)n} \sum_{v \in \cG} \sum_{1 \le r \le
   M}\sum_{\l:|k-\l| \le M} E\{J(v,k,r)I_n\}\\
&\sum_{0 \le k < (1-\eta)n}\sum_{\l:|k-\l| \le M} P\{\cE_n\} \le n(2M+1)u_n.
\end{split}
\]
These terms are therefore $o(n^2u_n)$. Similarly to \eqref{5.51a}
the remaining terms contribute
\begin{equation}
\begin{split}
&2 \sum_{0 \le k < (1-\eta)n} \sum_{v \in \cG} \sum_{1 \le r  \le
M \land \big((1-\eta)n - k\big)}
 \sum_{k+M < \l < (1-\eta)n} \sum_{w \in \cG}
\sum_{1 \le s  \le M \land \big((1-\eta)n - \l\big)}\\
&\phantom{MMMMMMMMMMMMMMMMMM}E\{J(v,k,r)J(w,\l,s)I_n\}\\
&= 2 \sum_{0 \le k < (1-\eta)n} \sum_{v \in \cG} \sum_{1 \le r
\le M \land \big((1-\eta)n - k\big)}
 \sum_{k+M < \l < (1-\eta)n} \sum_{w \in \cG}
\sum_{1 \le s  \le M \land \big((1-\eta)n - \l\big)}\\
&\phantom{MMMMMM}
P\{S_k = v\}f_rP\{S_\l = w|S_{k+r} = v\}f_sP\{S_n = e|S_{\l+s}=w\}\\
&=2 \sum_{0 \le k < (1-\eta)n} \sum_{1 \le r  \le M \land
\big((1-\eta)n - k\big)}\sum_{k+M < \l < (1-\eta)n} \sum_{1
  \le s \le M \land \big((1-\eta)n - \l\big)} f_rf_su_{n-r-s}.
\end{split}
\end{equation}
After division by $n^2u_n$ we find, just as in \eqref{1.19c},
\eqref{1.19d} that
\[
\limsup_{n \to \infty,p|n} \frac 1{n^2} E\{W^2_M|\cE_n\} \le
(1-\eta)^2\Big[\sum_{1 \le r \le M} f_r\rho^r\Big]^2.
\]
Together with \eqref{1.19f} this shows that
\begin{equation}
\lim_{n \to \infty, p|n} \frac 1{n^2}\text{ Var}\{W_M|\cE_n\} = 0,
\label{1.19h}
\end{equation}
and therefore, $(1/n)W_M$, conditioned on $\cE_n$, tends to $(1-\eta)
\sum_{r=1}^Mf_r \rho^r$ in probability as $n \to \infty$. By first
taking $M$ large and then $\eta$ small in
\eqref{1.19g} we obtain the desired \eqref{1.18b}. \qed

\medskip
\noindent {\it Example. Simple random walk on a regular tree}.
Here we shall explicitly calculate the values of $F$ and $F(\rho)$
which were stated in \eqref{1.18j}. We take an arbitrary vertex of
the tree $\cG_b$  for $e$. This will remain fixed throughout the
calculation. Also $\{S_n\}$ will be simple random walk on this
tree. Unless otherwise stated $S_0 = e$. This random walk has
period $p = 2$. For any vertex $v$, $d(v)$ denotes the number of
edges in the simple path on $\cG$ from $e$ to $v$; this is also
called the {\it height} of $v$. We set
\begin{equation}
T_n := d(S_n). \label{3.1}
\end{equation}
Again, unless stated otherwise $T_0 = d(S_0) = 0$.

It is well known and easy to prove that with $\{S_n\}$ simple
random walk on $\cG_b$, $\{T_n\}$ is a nearest neighbor random
walk on $\Bbb Z_+ = \{0,1,2,\dots\}$ with transition probabilities
\begin{equation}
P(x,y) = \begin{cases} \frac b{b+1} &\text{ if } x \ge 1, y = x+1\\
\frac 1{b+1} &\text{ if } x \ge 1, y = x-1\\
1 &\text{ if } x=0, y=1.
\end{cases}
\end{equation}
$P(x,y) = 0$ if $|x-y| > 1$ or $y < 0$. From this observation we
have
\begin{equation}
\begin{split}
f_r &= P\{\text{first return by $T_.$ to the origin is at time }r\},\\
u_r &= P\{T_r = 0\}.
\end{split}
\end{equation}
$f_r = u_r = 0$ if $r$ is odd, while explicit formulae for
$f_{2k}$ and for the generating functions of the $f_r$ and $u_r$
are known. Indeed, set
\[
\la =  \frac b{(b+1)^2}.
\]
Then, by the arguments in \cite{F} for the equations
XIII.4.6-XIII.4.8 (see also \cite {F}, Section XI.3)
\begin{equation}
F(z) := \sum_{r=1}^\infty f_rz^r = \frac {b+1}{2b} - \frac
{b+1}{2b} \big[1-4\la z^2\big]^{1/2},\; |z| \le 1,
\end{equation}
and hence, by recurrent event theory
\begin{equation}
\begin{split}
U(z) := \sum_{r=0}^\infty u_rz^r &= 2b \big[b-1 + (b+1)\big(1-4\la
z^2\big)^{1/2}\big]^{-1},\; |z| < 1.
\end{split}
\label{1.19j}
\end{equation}
Expansion of $F(z)$ shows that
\[
f_{2k} = \frac {b+1}{2b} \cdot \frac {(-1)^{k-1}}{2k-1} \binom{\frac 12}k
(4\la)^k = \frac{b+1}b \cdot \frac 1{2k-1} \binom{2k-1}k \la^k.
\]
Also,
\[
F = F(1) = \frac {b+1}{2b} - \frac{b+1}{2b} [1-4\la]^{1/2} = \frac
{b+1}{2b} - \frac{b+1}{2b}\cdot \frac{b-1}{b+1} = \frac 1b,
\]
as claimed in \eqref{1.18j}.

As pointed out in the beginning of the proof of Theorem 3,
$\rho = \min\{\rho_1, \rho_2\}$. In the present example $\rho_1$ is
the first place where $1-4\la x^2$ becomes 0, i.e.,
$\rho_1 = 1/(2 \sqrt \la)$. $\rho_2 > \rho_1$ because
$F(\rho_1) = F\big(1/(2\sqrt \la) = (b+1)/(2b)$, which is still
less than 1. Hence, $\rho =  1/(2\sqrt \la) = (b+1)/(2\sqrt b)$ and $F(\rho)
= (b+1)/(2b)$. This proves \eqref{1.18j}.
\qed

To conclude we prove the equivalence of \eqref{1.16a} and polynomial
growth of the group on which the random walk takes place.
\begin{lm} Let $(G, \cS)$ be a Cayley graph and let $\cV(n)$ be as in
\eqref{2.89ab}
Then, for a random walk on $(G, \cS)$ which satisfies \eqref{1.5} and
\eqref{1.6}, \eqref{1.16a} is equivalent to
\begin{equation}
\begin{split}
&\text{there exist some constants $C_i < \infty$ such that }\\
&\cV(n) \le C_{19}n^{C_{20}} \text{ for } n \ge 1.
\end{split}
\label{2.90}
\end{equation}
\end{lm}
\noindent
{\it Proof.}

Some version of this result was known to
N. Varopoulos. We learned the following argument from
Laurent Saloff-Coste.
Polynomial growth of $G$ as in \eqref{2.90}
implies \eqref{1.16a} by means of Theorem 5.1 in \cite{HeS}, as we
already observed for \eqref{1.18}.For the converse, assume
\eqref{1.16a} holds, and
let the constant $C_{21}$ be such that $u_{2n}/u_{4n} \le C_{21}$ for all
$n \ge 1$. It is shown in \cite{HeS}, equation (10) that
\begin{equation}
u_{2n+m} \le \frac 2{\cV\big(r(n,m)\big)}, n,m \ge 1
\label{2.91}
\end{equation}
where
\[
r(n,m) = \sqrt m \frac{u_{2n+m}}{u_{2n}}.
\]
If we take $m = 2n$ we get from \eqref{2.91}
\[
u_{4n} \le \frac 2{\cV\big(C_{21}\sqrt {2n}\big)}.
\]
But we already saw in \eqref{2.93} that \eqref{1.16a} implies $u_{2n}
\ge n^{-C_{22}}$ for some constant $C_{22}$ and large $n$. Thus
\[
\cV\big(C_{21} \sqrt{2n}\big) \le [2n]^{C_{22}} \text{ for large $n$},
\]
so that $\cV(n)$ cannot grow faster than a power of $n$.
\qed


\begin{thebibliography}{ABC}

\bibitem{Av} A. Avez, Limite de quotients pour les marches
  al\'eatoires sur les groupes, C. R. Acad. Sc. Paris, S\'er. A {\bf
  276} (1973) 317-320.

\bibitem{BJ} Ph. Bougerol and T. Jeulin, Brownian bridge on hyperbolic
  spaces and on homogeneous trees, Probab. Theory Relat. Fields {\bf 115}
  (1999) 95-120.

\bibitem{Ca} D. I. Cartwright, Some examples of random walks on free
  products of discrete groups, Ann. Mat. Pura Appl. {\bf 151} (1988) 1-15.

\bibitem{ChT} Y. S. Chow and H. Teicher, Probability Theory,
  $3^{\text{nd}}$ ed., Springer-Verlag, 1997.

\bibitem{De}
F. M. Dekking, On transience and recurrence of generalized random
walks, Z. Wahrsch. verw. Gebiete {\bf 61} (1982) 459-465.


\bibitem{DoV} M. D. Donsker and S. R. S.  Varadhan, On the number of distinct
  sites visited by a random walk, Comm. Pure Appl. Math. {\bf 32}
  (1979) 721-747.

\bibitem{DvE} A. Dvoretzky and P. Erd\H os, Some problems on random
  walk in space, Proc. Second Berkeley Symposium on Mathematical
Statistics and  Probability (J. Neyman ed.), pp. 353-368, Univ. of
California Press, 1951

\bibitem{E} A. Erschler, Isoperimetry for wreath products of Markov
  chains and multiplicity of selfintersections, Probab. Theory
  Relat. Fields, 2006.

\bibitem{F} W. Feller, An Introduction to Probability Theory and its
  Applications, vol I, $3^{\text{rd}}$ ed., John Wiley \& Sons, 1968.

\bibitem{Ge} P. Gerl, Diskrete, mittelbare Gruppen,
  Monatsh. Math. {\bf 77} (1973) 307-318.


\bibitem{Gr} R. I. Grigorchuk, Degrees of growth of finitely generated
groups and the theory of invariant means, Math. USSR Izvestiya {\bf
  25} (1985) 259-300.

\bibitem{Ha} Y. Hamana, An almost sure invariance principle
  for the range of random walks, Stochastic Process. Appl. {\bf 78}
  (1998) 131-143,

\bibitem{HeS} W. Hebisch and L. Saloff-Coste, Gaussian estimates for
  Markov chains and random walks on groups, Ann. Probab. {\bf 21}
  (1993) 673-709.

\bibitem{JP} N. C. Jain and W. E. Pruitt, The range of random walk, Proc.
Sixth Berkeley Symposium on Mathematical Statistics and Probability
(L. M. Le Cam, J. Neyman and E. L. Scott eds.) Vol 3, pp. 31-50, Univ. of
California Press, 1972.

\bibitem{KV} V. A. Kaimanovich and A. M. Vershik, Random walks on
  discrete groups: boundary and entropy, Ann. Probab. {\bf 11} (1983)
  457-490.

\bibitem{K1} H. Kesten, Full Banach mean values on countable groups,
  Math. Scan. {\bf 7} (1959) 146-156.

\bibitem{K2} H. Kesten, Subdiffusive behavior of random
walk on a random cluster, Ann. Inst. Henri Poincar\'e, Probabilit\'es
et Statistiques {\bf 22} (1986) 425-487.

\bibitem{L} S. P. Lalley, Finite range random walk on free groups and
  homogeneous trees, Ann. Probab. {\bf 21} (1993) 2087-2130.

\bibitem{LP} R. Lyons and Y. Peres, Probability on Trees and
  Networks. To appear
(2006); see http://mypage.iu.edu/~rdlyons/prbtree/prbtree.html.

\bibitem{NGV} S. K. Nechaev, A. Yu. Grosberg and A. M. Vershik, Random
  walks on braid groups: Brownian bridges, complexity and statistics,
  J. Physics A Math. Gen. {\bf 29} (1996) 2411-2434.

\bibitem{Pa} A. L. T. Paterson, Amenability, Math Surveys and Monographs
  {\bf 29}, Amer. Math. Soc. 1988.

\bibitem{PS} C. Pittet and L. Saloff-Coste, On random walks on wreath
  products, Ann. Probab. {\bf 30} (2002) 948-977.

\bibitem{R} W. Rudin, Functional Analysis, $2^{\text{nd}}$ ed.,
  Mc-Graw Hill, 1991.

\bibitem{Sa} S. Sawyer, Isotropic random walks in a tree,
  Z. Wahrsch. verw. Gebiete {\bf 42} (1978) 279-292.

\bibitem{Si} A. Sinclair,
Algorithms for Random Generation and Counting: a Markov Chain
Approach, Birkhauser Verlag, Boston, 1993.

\bibitem{Sp1} F. Spitzer, Principles of Random Walk, $2^{\text{nd}}$ ed.,
Springer-Verlag, 1976.

\bibitem{Sp2} F. Spitzer in Discussion of Subadditive ergodic theory by
  J. F. C. Kingman, Ann Probab. {\bf 1} (1973) 904-905.


\bibitem{V}
D. Vere-Jones, Ergodic properties of nonnegative matrices I,
Pacific J. Math. {\bf 22} (1967) 361-386.

\bibitem{Wo}
W. Woess,  Random Walks on Infinite Graphs and Groups. Cambridge
Univ. Press, 2000.

\bibitem{Y}
R. Young, Averaged Dehn functions for nilpotent groups, arXiv math.
GR/0510665.

\end{thebibliography}
\end{document}